# General Eulerian Numbers and Eulerian Polynomials

Tingyao Xiong[1]
txiong@radford.edu

Hung-ping Tsao[2]
tsaohp.tsao6@gmail.com

Jonathan I. Hall[3]
jhall@math.msu.edu

*Abstract*  In this paper, we will define general Eulerian numbers and Eulerian polynomials based on general arithmetic progressions. Under the new definitions, we have been successful in extending several well-known properties of traditional Eulerian numbers and polynomials to the general Eulerian polynomials and numbers.

*Keywords*: Eulerian numbers, Eulerian polynomials, generating functions.

## 1. Introduction

Jacques Bernoulli ([1], page 95-97) had introduced his famous *Bernoulli numbers*, denoted by $B_{2r}$ ($B_{2r+1} = 0$ for $r \geq 1$) to evaluate the sum of the $n$-th power of the first $m$ integers. He then proved the following summation formula

$$\sum_{i=1}^{m} i^n = \frac{m^{n+1}}{n+1} + \frac{m^n}{2} + \frac{1}{n+1}\sum_{r=1}^{\lfloor n/2 \rfloor} \binom{n+1}{2r} m^{n-2r+1}(-1)^{r+1} B_{2r} \qquad (1)$$

when $n, m \geq 1$.

Two decades later, Euler ([2]) studied the *alternating sum* $\sum_{i=1}^{m}(-1)^i i^n$. He ended up with giving the following general result ([8], (2.8), page 259)

$$\sum_{i=1}^{m} i^n t^i = \sum_{l=1}^{n}(-1)^{n+l}\binom{n}{l}\frac{t^{m+1}m^l}{(t-1)^{n-l+1}}A_{n-l}(t) + (-1)^n \frac{t(t^m - 1)}{(t-1)^{n+1}} A_n(t) \qquad (2)$$

Another simplified form of $\sum_{i=1}^{m} i^n t^i$ is the following ([8], (3.3), page 263) :

$$\sum_{i=1}^{m} i^n t^i = -t^{m+1}\sum_{k=0}^{n}\binom{n}{k}\frac{m^{n-k}}{(1-t)^{k+1}}A_k(t) + \frac{t}{(1-t)^{n+1}}A_n(t) \qquad (3)$$

where $A_n(t), \ (n = 0, 1, 2, \ldots)$ are called *Eulerian Polynomials*, and are recursively defined by

$$A_0(t) = 1 \quad \text{and} \quad A_n(t) = \sum_{k=0}^{n-1} \binom{n}{k} A_k(t)(t-1)^{n-1-k} \tag{4}$$

---
[1] Radford University
[2] Novato, CA
[3] Michigan State University

([8], (2.7), page 264)

Note that for $\sum_{i=1}^{m} i^n t^i$, if we put $t = -1$, then $\sum_{i=1}^{m} i^n t^i$ becomes the alternating sum of the $n$-th power of the first $m$ positive integers. Furthermore, as $m \to \infty$, we have ([8], (3.2), page 263)

$$\frac{A_n(t)}{(1-t)^{n+1}} = \sum_{j=0}^{\infty} t^j (j+1)^n \quad (n \geq 0) \tag{5}$$

Each Eulerian polynomial can be presented as a generating function of *Eulerian numbers* $A_{n,k}$ ([3]), also introduced by Euler, as

$$A_n(t) = \sum_{k=0}^{n-1} A_{n,k} t^k \tag{6}$$

Furthermore, the corresponding exponential generating function ([8], (3.1), page 262) is

$$A(t, u) = \sum_{n \geq 0} A_n(t) \frac{u^n}{n!} = \frac{t-1}{t - \exp(u(t-1))} \tag{7}$$

The following combinatorial definition of Eulerian numbers was discovered by Riordan in the 1950's:

**Definition 1.1** *A Eulerian number $A_{n,k}$ is the number of permutations $p_1 p_2 p_3 \cdots p_n$ of the first $n$ numbers $\{1, 2, \ldots, n\}$ that have $k$ ascents (or descents), i.e., $k$ places where $p_j < p_{j+1}$ (or $p_j > p_{j+1}$), $1 \leq j \leq n-1$.*

**Example 1** *When $n = 3$, $k = 1$, $A_{3,1} = 4$, because there are four permutations with only one ascent:*

$$132, \quad 213, \quad 231, \quad 312$$

Then $A_{n,k}$ satisfies the recurrence: $A_{n,0} = 1$, $(n \geq 0)$, $A_{n,k} = 0$ $(k \geq n)$ and

$$A_{n,k} = (k+1)A_{n-1,k} + (n-k)A_{n-1,k-1} \quad (1 \leq k \leq n-1) \tag{8}$$

It is well known that Eulerian numbers have the following symmetric property:

**Proposition 1.2** *Given a positive integer $n$, and $0 \leq k \leq n-1$, $A_{n,k} = A_{n,n-1-k}$.*

**Proof.** It is obvious because if a permutation has $k$ ascents then its reverse has $n-1-k$ ascents. □

Besides the recursive formula (8), $A_{n,k}$ can be calculated directly by the following analytic formula ([8], (3.5), page 264):

$$A_{n,k} = \sum_{i=0}^{k} (-1)^i (k-i+1)^n \binom{n+1}{i} \quad (0 \leq k \leq n-1) \tag{9}$$

Since the 1950's, Carlitz ([6], [7]) and his successors have generalized Euler's results to $q$-sequences $\{1, q, q^2, q^3, \ldots\}$. Under Carlitz's definition, the $q$-Eulerian numbers $A_{n,k}(q)$ are given by

$$[x]^n = \sum_{k=0}^{n-1} A_{n,k}(q) \begin{bmatrix} x+k-1 \\ n \end{bmatrix} \quad (n \geq 1), \tag{10}$$

where

$$[x] = \frac{1-q^x}{1-q} \tag{11}$$

and

$$\begin{bmatrix} x \\ n \end{bmatrix} = \frac{(1-q^x)(1-q^{x-1})\ldots(1-q^{x-n+1})}{(1-q)(1-q^2)\ldots(1-q^n)} \tag{12}$$

Then the $q$-Eulerian polynomials $A_n(t,q)$ are defined as

$$A_n(t,q) = \sum_{k=0}^{n-1} A_{n,k} t^k \quad (n \geq 1) \tag{13}$$

Like the traditional Eulerian numbers, the $q$-Eulerian numbers $A_{n,k}(q)$ have the following recursive formula

$$\begin{aligned} & A_{1,0}(q) = 1, A_{n,k}(q) = 0 \quad \text{if} \quad (k \leq 0 \quad \text{or} \quad k \geq n). \\ & A_{n,k}(q) = q^{n-1-k}[k+1]A_{n-1,k}(q) + [n-k]A_{n-1,k-1}(q) \quad (n \geq 1) \end{aligned} \tag{14}$$

**Definition 1.3** *Given a permutation* $\pi = p_1 p_2 p_3 \cdots p_n$ *of the first* $n$ *positive integers, define*

$$\text{maj } \pi = \sum_{p_j > p_{j+1}} j \quad \text{and}$$

$$a(n,k,i) = \#\{\pi \mid \text{maj } \pi = i \ \& \ \pi \text{ has } k \text{ ascents}\}$$

In 1974, Carlitz ([7]) completed his study of his $q$-Eulerian numbers by giving a combinatorial meaning to his $q$-Eulerian numbers:

$$A_{n,k}(q) = q^{\frac{(n-k+1)(n-k)}{2}} \sum_{i=0}^{k(n-k-1)} a(n, n-k, i) q^i \tag{15}$$

where functions $a(n,k,i)$ are as defined in Definition 1.3. Interested readers can find more detail about the history of Eulerian numbers, Eulerian polynomials and the corresponding concepts in $q$-environment from [8].

In this paper, instead of studying $q$-sequences, we will generalize Euler's work on Eulerian numbers and Eulerian polynomials to any general arithmetic progression

$$\{a, a+d, a+2d, a+3d, \ldots\}$$

In section 2, we will give a new definition of general Eulerian numbers based on a given arithmetic progression as defined in (16). Under the new definition, some well-known combinatorial properties of traditional Eulerian numbers become special cases of our more general results. In section 3, we will define general Eulerian polynomials. Then equations (2) –(9) become special cases of our more general results.

## 2. General Eulerian Numbers

The traditional Eulerian numbers $A_{n,k}$, play an important role in the well-known Worpitzky's Identity ([9]):

$$x^n = \sum_{k=0}^{n-1} \binom{x+k}{n} A_{n,k} \quad n \geq 1 \tag{17}$$

Before we give a general definition of Eulerian numbers based on a given arithmetic progression (16), we shall mention a property associated with the traditional Eulerian numbers $A_{n,k}$.

**Proposition 2.1** Let $A_{n,k}$ be as defined in Definition 1.1, then

$$\sum_{i=1}^{m} i^n = \sum_{k=0}^{n-1} A_{n,k} \binom{m+k+1}{n+1}$$

**Proof.** see [10], (4), page 348. □

Given an arithmetic progression (16), we want to define general Eulerian numbers $A_{n,k}(a,d)$ so that the important properties of traditional Eulerian numbers such as the recursive formula (8), Worpitzky's Identity (17), and Proposition 2.1, etc, become special cases of more general results under the new definition.

**Definition 2.2** *General Eulerian numbers $A_{n,k}(a,d)$ associated with an arithmetic progression as in (16) are defined as* $A_{0,-1} = 1$, $A_{n,k} = 0$ $(k \geq n \text{ or } k \leq -2)$ *and*

$$A_{n,k}(a,d) = (-a+(k+2)d)A_{n-1,k}(a,d) + (a+(n-k-1)d)A_{n-1,k-1}(a,d) \quad (1 \leq k \leq n-1) \tag{18}$$

Note that the traditional Eulerian numbers have become the special cases of general eulerian numbers when $a = d = 1$ in the arithmetic progression

$$\{a, a+d, a+2d, a+3d, \ldots\}$$

In other words, the traditional Eulerian numbers correspond to the sequence of natural numbers

$$\{1, 2, 3, 4, \ldots\}$$

With the new Definition 2.2, we are able to prove the following two properties.

**Lemma 2.3 (General Worpitzky's Identity)** Given an arithmetic progression as in (16),

$$(a + (i-1)d)^n = \sum_{j=-1}^{n-1} A_{n,j}(a,d) \binom{i+j}{n} \tag{19}$$

**Proof.** We will prove Lemma 2.3 by induction on $n$.

When $n = 1$,

$$a + (i-1)d = A_{1,-1}(a,d)\binom{i-1}{1} + A_{1,0}(a,d)\binom{i}{1} = a + (i-1)d$$

Now suppose

$$(a + (i-1)d)^n = \sum_{j=-1}^{n-1} A_{n,j}(a,d) \binom{i+j}{n} \quad \text{then } (a + (i-1)d)^{n+1}$$

$$= \left[\sum_{j=-1}^{n-1} A_{n,j}(a,d) \binom{i+j}{n}\right](a + (i-1)d)$$

$$= \sum_{j=-1}^{n-1} A_{n,j}(a,d)\left[a\binom{i+j}{n} - (j+2)d\binom{i+j}{n} + (i+j+1)d\binom{i+j}{n}\right]$$

$$= \sum_{j=-1}^{n-1} A_{n,j}(a,d)\left[a\binom{i+j+1}{n+1} - a\binom{i+j}{n+1} + (n-j-1)d\binom{i+j+1}{n+1}\right]$$

$$+ \sum_{j=-1}^{n-1} A_{n,j}(a,d)\left[(j+2)d\binom{i+j}{n+1}\right]$$

$$= \sum_{j=-1}^{n-1} A_{n,j}(a,d)(a + (n-j-1)d)\binom{i+j+1}{n+1} + \sum_{j=-1}^{n-1} A_{n,j}(a,d)(-a + (j+2)d)\binom{i+j}{n+1}$$

$$= \sum_{j=-1}^{n} \left[A_{n,j-1}(a,d)(a + (n-j)d) + A_{n,j}(a,d)(-a + (j+2)d)\right]\binom{i+j}{n+1}$$

$$= \sum_{j=-1}^{n} A_{n+1,j}(a,d)\binom{i+j}{n+1}$$

□

With Lemma 2.3, the following Lemma which is a generalization of Proposition 2.1.

**Lemma 2.4** Given an arithmetic progression as in (16),

$$\sum_{i=1}^{m}(a + (i-1)d)^n = \sum_{j=-1}^{n-1} A_{n,j}(a,d)\binom{m+j+1}{n+1} \tag{20}$$

**Proof.** We will prove Lemma 2.4 by induction on $m$.

When $m = 1$,

$$a + (1-1)d = A_{1,0}(a,d)\binom{1}{1} = a \times 1 = a$$

Now suppose

$$\sum_{i=1}^{m}(a + (i-1)d)^n = \sum_{j=1}^{n-1} A_{n,j}(a,d)\binom{m+j+1}{n+1}$$

Then from Definition 2.2 and Lemma 2.3

$$\sum_{i=1}^{m+1}(a + (i-1)d)^n = \sum_{j=1}^{n-1} A_{n,j}(a,d)\binom{m+j+1}{n+1} + (a+md)^n$$

$$= \sum_{j=1}^{n-1} A_{n,j}(a,d)\left[\binom{m+j+1}{n+1} + \binom{m+j+1}{n}\right]$$

$$= \sum_{j=1}^{n-1} A_{n,j}(a,d)\binom{m+1+j+1}{n+1}$$

□

The general Eulerian numbers $A_{n,k}(a,d)$ can be calculated directly from the following formula, which is a generalization of equation (9).

**Lemma 2.5** For a given arithmetic progression $\{a, a+d, a+2d, a+3d, \ldots\}$, the general Eulerian numbers satisfy

$$A_{n,k}(a,d) = \sum_{i=0}^{k+1}(-1)^i[(k+2-i)d - a]^n \binom{n+1}{i}$$

**Proof.** Again, the proof is by induction on $n$.

For $n = 0, k = -1$,
$$A_{0,-1}(a,d) = (-1)^0(d-a)^0 \binom{1}{0} = 1$$

For $n = 1, k = -1$,
$$A_{1,-1}(a,d) = (-1)^0(d-a)^1\binom{2}{0} = d - a$$

For $n = 1, k = 0$,
$$A_{1,0}(a,d) = (-1)^0(2d-a)^1\binom{2}{0} + (-1)^1(d-a)^1\binom{2}{1} = a$$

Now suppose

$$A_{n-1,k}(a,d) = \sum_{i=0}^{k+1}(-1)^i[(k+2-i)d-a]^{n-1}\binom{n}{i},$$

Then from the recursive formula (18),

$$A_{n,k}(a,d) = [-a+(k+2)d]A_{n-1,k}(a,d) + [a+(n-k-1)d]A_{n-1,k-1}(a,d)$$

$$= [-a+(k+2)d]\sum_{i=0}^{k+1}(-1)^i[(k+2-i)d-a]^{n-1}\binom{n}{i}$$

$$+ [a+(n-k-1)d]\sum_{i=1}^{k+1}(-1)^{i-1}[(k+2-i)d-a]^{n-1}\binom{n}{i-1} \quad \text{by induction}$$

$$= [(k+2)d-a]\binom{n}{0} + \sum_{i=1}^{k+1}(-1)^i[(k+2-i)d-a]^{n-1}$$

$$\left[(k+2)d\binom{n}{i} - (n-k-1)d\binom{n}{i-1} - a\binom{n}{i} - a\binom{n}{i-1}\right]$$

$$= \sum_{i=0}^{k+1}(-1)^i[(k+2-i)d-a]^n\binom{n+1}{i}$$

□

## 3. General Eulerian Polynomials

Definition 3.1 *We define the general Eulerian Polynomials associated to an arithmetic progression as in (16) as*

$$T_n(t,a,d) = \sum_{k=-1}^{n-1} A_{n,k}(a,d)\, t^{k+1}$$

Definition 3.1 is a generalization of the traditional Eulerian Polynomials as in (6). The following lemma gives the relation between the general Eulerian Polynomials and the traditional Eulerian Polynomials.

Lemma 3.2  Let $T_n(t,a,d)$ be the general Eulerian Polynomials as in Definition 3.1, $T_0(t,a,d) = 1$. Then

$$T_n(t,a,d) = \sum_{k=-1}^{n-1} A_{n,k}(a,d)\, t^{k+1} = \sum_{j=0}^{n}\binom{n}{j} d^j A_j(t)(at-a)^{n-j}$$

where $A_j(t)$, $j = 0, 1, 2, \ldots, n$ are traditional Eulerian polynomials as defined in equations (4) and (6).

**Proof.** When $n = 0$,

$$T_0(t,a,d)\begin{pmatrix}0\\0\end{pmatrix}d^0A_0(t)(at-a)^0 = 1$$

Now suppose

$$T_n(t,a,d) = \sum_{j=0}^{n}\begin{pmatrix}n\\j\end{pmatrix}d^jA_j(t)(at-a)^{n-j}$$ Then

$$T_{n+1}(t,a,d) = \sum_{k=1}^{n} A_{n+1,k}(a,d)t^{k+1} \quad \text{by definition}$$

$$= \sum_{k=1}^{n}(a+(n-k)d)A_{n,k-1}(a,d)t^{k+1} + \sum_{k=1}^{n}(-a+(k+2)d)A_{n,k}(a,d)t^{k+1}$$

$$= \sum_{k=1}^{n-1}(a+(n-k-1)d)A_{n,k}(a,d)t^{k+2} + \sum_{k=1}^{n-1}(-a+(k+2)d)A_{n,k}(a,d)t^{k+1}$$

$$= \sum_{k=1}^{n-1}(at-a)A_{n,k}(a,d)t^{k+1} + \sum_{k=1}^{n-1}d(nt-t+2)A_{n,k}(a,d)t^{k+1}$$

$$+ \sum_{k=1}^{n-1}d(1-t)kA_{n,k}(a,d)t^{k+1}$$

$$= (at-a)T_n(t,a,d) + (dnt-dt+2d)T_n(t,a,d) + \sum_{k=1}^{n-1}d(1-t)kA_{n,k}(a,d)t^{k+1}$$

(21)

Note that by definition, $T_n(t,a,d) = \sum_{k=1}^{n-1}A_{n,k}(a,d)t^{k+1}$. So
$$T'_n(t,a,d) = \sum_{k=1}^{n-1}(k+1)A_{n,k}(a,d)t^k,$$
which implies

$$tT'_n(t,a,d) = T_n(t,a,d) + \sum_{k=1}^{n-1}kA_{n,k}t^{k+1} \tag{22}$$

On the other hand, from ([8], (3.4), page 263), we have

$$A_n(t) = [1+(n-1)t]A_{n-1}(t) + t(1-t)A'_{n-1}(t) \tag{23}$$

With results as in (22) and (23), we have

$$(dnt - dt + 2d)T_n(t,a,d) + \sum_{k=1}^{n-1} d(1-t)kA_{n,k}(a,d)t^{k+1}$$

$$= (dnt + d)T_n(t,a,d) + d(1-t)tT'_n(t,a,d)$$

$$= (dnt + d)T_n(t,a,d) + d(1-t)t\sum_{j=0}^{n} d^j \binom{n}{j} \times$$

$$[A_j^{'n-j} + a(n-j)A_j(t)(at-a)^{n-j-1}] \qquad \text{by induction}$$

$$= (dnt + d)T_n(t,a,d) - dt(n-j)T_n(t,a,d) + d(1-t)t \sum_{j=0}^{n} \binom{n}{j} d^j A_j^{'n-j}$$

$$= \sum_{j=0}^{n} \binom{n}{j} d^{j+1}(at-a)^{n-j}[(1+jt)A_j(t) + t(1-t)A'_j(t)]$$

$$= \sum_{j=0}^{n} \binom{n}{j} d^{j+1}(at-a)^{n-j}A_{j+1}(t)$$

Therefore, expression (21) becomes

$$T_{n+1}(t,a,d) = \sum_{j=0}^{n} \binom{n}{j} d^{j+1}(at-a)^{n-j}A_{j+1}(t) + \sum_{j=0}^{n} \binom{n}{j} d^j(at-a)^{n-j+1}A_j(t)$$

$$= \sum_{j=1}^{n+1} \binom{n}{j-1} d^j(at-a)^{n+1-j}A_j(t) + \sum_{j=0}^{n} \binom{n}{j} d^j(at-a)^{n-j+1}A_j(t)$$

$$= \sum_{j=1}^{n} \left[\binom{n}{j} + \binom{n}{j+1}\right] d^j(at-a)^{n+1-j}A_j(t)$$

$$+ \binom{n}{n} d^{n+1}(at-a)^0 A_{n+1}(t) + \binom{n}{0} d^0(at-a)^{n+1}A_0(t)$$

$$= \sum_{j=0}^{n+1} \binom{n+1}{j} d^j A_j(t)(at-a)^{n+1-j}$$

□

The following result is a generalization of equation (7).

**Lemma 3.3** Let $T_n(t,a,d)$ be as defined in Lemma 3.2. Then

$$T_n(t,u,a,d) = \sum_{n \geq 0} T_n(t,a,d)\frac{u^n}{n!} = \frac{(t-1)\exp(au(t-1))}{t - \exp(du(t-1))}$$

**Proof.** From Lemma 3.2,

$$T_n(t,u,a,d) = \sum_{n\geq 0} T_n(t,a,d)\frac{u^n}{n!} = \sum_{n\geq 0}\left[\sum_{j=0}^{n}\binom{n}{j}d^j A_j(t)(at-a)^{n-j}\right]\frac{u^n}{n!}$$

$$= \sum_{n\geq 0}\left[\sum_{\substack{k+l=n \\ 0\leq k\leq n}}\frac{n!}{k!\,l!}d^k A_k(t,a,d)(at-a)^l\frac{u^k \cdot u^l}{n!}\right] = \left(\sum_{k\geq 0}A_k(t)\frac{(du)^k}{k!}\right)\left(\sum_{l\geq 0}\frac{(ua(t-1))^l}{l!}\right)$$

$$= \frac{(t-1)\exp(au(t-1))}{t-\exp(du(t-1))}$$

□

Using the results from Lemma 3.3, we can derive the following lemma, which is a general version of equation (5).

**Proposition 3.4** Given an arithmetic progression as in (16), let $T_n(t, a-d, -d)$ be the general Eulerian Polynomials associated to the arithmetic progression $\{a-d,\ a-2d,\ a-3d,\ \ldots\}$. Then

$$-\frac{T_n(t, a-d, -d)}{(t-1)^{n+1}} = \sum_{j=0}^{\infty} t^j (a+jd)^n \qquad (n \geq 0)$$

**Proof.**

$$\frac{e^{au}}{1-te^{du}} = e^{au}\sum_{j=0}^{\infty}(te^{du})^j = \sum_{j=0}^{\infty} t^j e^{u(a+jd)}$$

$$= \sum_{j=0}^{\infty} t^j \sum_{n=0}^{\infty}\frac{u^n(a+jd)^n}{n!} = \sum_{n=0}^{\infty}\frac{u^n}{n!}\left(\sum_{j=0}^{\infty} t^j (a+jd)^n\right)$$

On the other hand, by Lemma 3.3,

$$\sum_{n=0}^{\infty} -\frac{T_n(t, a-d, -d)}{(t-1)^{n+1}}\frac{u^n}{n!} = \frac{1}{t-1}\sum_{n=0}^{\infty} -\frac{T_n(t, a-d, -d)}{(t-1)^n}\frac{u^n}{n!}$$

$$= -\frac{1}{t-1}\sum_{n=0}^{\infty}\frac{T_n(t, a-d, -d)}{n!}\left(\frac{u}{t-1}\right)^n = -\frac{1}{t-1}\cdot\frac{(t-1)e^{(a-d)u}}{t-e^{-du}} = -\frac{e^{(a-d)u}}{t-e^{-du}} = \frac{e^{au}}{1-te^{du}}$$

By comparing the coefficients of $\frac{u^n}{n!}$, we have

$$-\frac{T_n(t, a-d, -d)}{(t-1)^{n+1}} = \sum_{j=0}^{\infty} t^j (a+jd)^n \qquad (n \geq 0)$$

□

For the finite summation $\sum_{i=1}^{m} t^i[a+(i-1)d]^n$, we have the following property which is a generalization of equations (2) and (3).

**Lemma 3.5** Let $T_n$ be the general Eulerian polynomials as defined in Definition 3.1

$$\sum_{i=1}^{m} t^i[a+(i-1)d]^n = \sum_{l=0}^{n} \binom{n}{l} \frac{t^{m+1}(dm-d)^l}{(t-1)^{n-l+1}} T_{n-l}(t,a,-d) - \frac{t^2}{(t-1)^{n+1}} T_n(t,a,-d) \quad (24)$$

and

$$\sum_{i=1}^{m} t^i[a+(i-1)d]^n = \frac{t^{m+1}}{(t-1)^{n+1}} T_n(t, a+d(m-1),-d) - \frac{t^2}{(t-1)^{n+1}} T_n(t,a,-d) \quad (25)$$

**Proof.** We will use equations (2) and (3) to prove expressions (24) and (25) respectively.

$$\sum_{i=1}^{m} t^i[a+(i-1)d]^n = \sum_{i=1}^{m} t^i \left[ \sum_{j=0}^{n} \binom{n}{j} a^{n-j}(i-1)^j d^j \right]$$

$$= \sum_{j=0}^{n} \binom{n}{j} d^j a^{n-j} \left[ \sum_{i=1}^{m} t^i (i-1)^j \right] = \sum_{j=0}^{n} \binom{n}{j} d^j a^{n-j} t \left[ \sum_{i=1}^{m-1} i^j t^i \right]$$

For equation (24), if we use equation (2) to evaluate $\sum_{i=1}^{m-1} i^j t^i$, we have

$$\sum_{i=1}^{m} t^i[a+(i-1)d]^n$$

$$= \sum_{j=0}^{n} \binom{n}{j} d^j a^{n-j} t \left[ \sum_{l=1}^{j} (-1)^{j+l} \binom{j}{l} \frac{t^m(m-1)^l}{(t-1)^{j-l+1}} A_{j-l}(t) + (-1)^j \frac{t^m - t}{(t-1)^{j+1}} A_j(t) \right]$$

$$= \sum_{j=0}^{n} \binom{n}{j} d^j a^{n-j} t \left[ \sum_{l=0}^{j} (-1)^{j+l} \binom{j}{l} \frac{t^m(m-1)^l}{(t-1)^{j-l+1}} A_{j-l}(t) \right]$$

$$+ \sum_{j=0}^{n} \binom{n}{j} d^j a^{n-j} t (-1)^j \frac{-t}{(t-1)^{j+1}} A_j(t) = I + II$$

$$II = \sum_{j=0}^{n} \binom{n}{j} d^j a^{n-j} t (-1)^j \frac{-t}{(t-1)^{j+1}} A_j(t)$$

$$= \frac{-t^2}{(t-1)^{n+1}} \sum_{j=0}^{n} \binom{n}{j} (-d)^j (a(t-1))^{n-j} A_j(t) = -\frac{t^2}{(t-1)^{n+1}} T_n(t,a,-d)$$

This gives the second term in equation (24).

$$I = \sum_{j=0}^{n} \binom{n}{j} d^j a^{n-j} t \sum_{l=0}^{j} (-1)^{j+l} \binom{j}{l} \frac{t^m(m-1)^l}{(t-1)^{j-l+1}} A_{j-l}(t)$$

$$= \sum_{l=0}^{n} \sum_{j=l}^{n} \binom{n}{j} d^j a^{n-j} t (-1)^{j+l} \binom{j}{l} \frac{t^m(m-1)^l}{(t-1)^{j-l+1}} A_{j-l}(t)$$

$$= \sum_{l=0}^{n} \sum_{j=l}^{n} \frac{n!}{j!(n-j)!} d^j a^{n-j} t (-1)^{j+l} \frac{j!}{l!(j-l)!} \frac{t^m(m-1)^l}{(t-1)^{j-l+1}} A_{j-l}(t)$$

$$= \sum_{l=0}^{n} \sum_{k=0}^{n-l} \frac{n!}{(n-k-l)!} d^{k+l} a^{n-k-l} t (-1)^k \frac{1}{l!k!} \frac{t^m(m-1)^l}{(t-1)^{k+1}} A_k(t) \qquad k = j-l$$

$$= \sum_{l=0}^{n} \frac{(m-1)^l n! t^{m+1} d^l}{l!(n-l)!} \sum_{k=0}^{n-l} \frac{(-d)^k a^{n-k-l}(n-l)!}{k!(n-k-l)!(t-1)^{k+1}} A_k(t)$$

$$= \sum_{l=0}^{n} \binom{n}{l} \frac{(m-1)^l t^{m+1} d^l}{(t-1)^{n-l+1}} \sum_{k=0}^{n-l} \binom{n-l}{k} (-d)^k (at-a)^{n-k-l} A_k(t)$$

$$= \sum_{l=0}^{n} \binom{n}{l} \frac{t^{m+1}(dm-d)^l}{(t-1)^{n-l+1}} T_{n-l}(t, a, -d)$$

which gives the first term of equation (24). So we have proved equation (24) by expression (2).

For equation (25), if we use equation (3) to evaluate $\sum_{i=1}^{m-1} i^j t^i$, we have

$$\sum_{i=1}^{m} t^i [a + (i-1)d]^n$$

$$= \sum_{j=0}^{n} \binom{n}{j} d^j a^{n-j} t \left[ -t^m \sum_{k=0}^{j} \binom{j}{k} \frac{(m-1)^{j-k}}{(1-t)^{k+1}} A_k(t) + \frac{t}{(1-t)^{j+1}} A_j(t) \right]$$

$$= -t^{m+1} \sum_{j=0}^{n} \binom{n}{j} d^j a^{n-j} \sum_{k=0}^{j} \binom{j}{k} \frac{(m-1)^{j-k}}{(1-t)^{k+1}} A_k(t) + \sum_{j=0}^{n} \binom{n}{j} d^j a^{n-j} (1-t)^{-j} \frac{t^2}{1-t} A_j(t)$$

$$= III + IV$$

$$IV = \sum_{j=0}^{n} \binom{n}{j} d^j a^{n-j} (1-t)^{-j} \frac{t^2}{1-t} A_j(t)$$

$$= \sum_{j=0}^{n} \binom{n}{j} (-d)^j a^{n-j} (t-1)^{-j} \frac{t^2}{(1-t)(t-1)^n} A_j(t) = -\frac{t^2}{(t-1)^{n+1}} T_n(t,a,-d)$$

by Lemma 3.2. Thus we have obtained the second part of equation (25).

$$III = -t^{m+1} \sum_{j=0}^{n} \binom{n}{j} d^j a^{n-j} \sum_{k=0}^{j} \binom{j}{k} \frac{(m-1)^{j-k}}{(1-t)^{k+1}} A_k(t)$$

$$= -t^{m+1} \sum_{j=0}^{n} \sum_{k=0}^{j} \frac{n!}{j!(n-j)!} \cdot \frac{j!}{k!(j-k)!} d^j a^{n-j} \frac{(m-1)^{j-k}}{(1-t)^{k+1}} A_k(t)$$

$$= -t^{m+1} \sum_{k=0}^{n} \frac{n!}{k!(1-t)^{k+1}} A_k(t) \sum_{j=k}^{n} \frac{1}{(n-j)!(j-k)!} (m-1)^{j-k} d^j a^{n-j}$$

$$= -t^{m+1} \sum_{k=0}^{n} \frac{n!}{k!(1-t)^{k+1}} A_k(t) \sum_{l=0}^{n-k} \frac{1}{l!(n-k-l)!} (m-1)^l d^{l+k} a^{n-l-k} \qquad l = j-k$$

$$= -t^{m+1} \sum_{k=0}^{n} \binom{n}{k} \frac{d^k}{(1-t)^{k+1}} A_k(t) (a + d(m-1))^{n-k}$$

$$= \frac{t^{m+1}}{(t-1)^{n+1}} \sum_{k=0}^{n} \binom{n}{k} (-d)^k (t-1)^{n-k} A_k(t) (a + d(m-1))^{n-k}$$

$$= \frac{t^{m+1}}{(t-1)^{n+1}} T_n(t, a + d(m-1), -d)$$

which gives the first term of equation (25). □

The authors want to thank Professor Dominique Foata for his kind help.

scientiarum imperialis Petropolitana 8, page 147-158, 1736.

[4] J. Riordan, "An Introduction to Combinatorial Analysis", *New York, Wiley*, 1958.

[5] Miklós Bóna, "Combinatorics of Permutations", *Chapman & Hall/CRC*, 2004.

[6] L. Carlitz, "$q-$Bernoulli and Eulerian numbers", *Transaction of American Mathematical Society*, Vol. 76, page 332--350, 1954.

[7] L. Carlitz, "A Combinatorial Property of $q-$Eulerian Numbers", *American Mathematical Monthly*, Vol. 82, page 51--54, 1975.

[8] D. Foata, "Eulerian Polynomial: From Euler's Time to the Present", *The Legacy of Alladi Ramakrishnan in the Mathematical Sciences*, Springer Science+Business Media, page 253-273, LLC 2010.

[9] J. Worpitzky, "Studien über die Bernoullischen und Eulerischen Zahle", *J. reine angew. Math*, Vol. 94, page 203-232, 1883.

[10] H. Tsao, "Sums of Powers and Eulerian Number", *Mathematical Gasette*, Vol. 95, page 347-349, 2011.
**Address:**
Department of Mathematics and Statistics
Radford University
Radford, VA 24141, U.S.A.